\documentclass{article}
\usepackage[leqno]{amsmath}
\usepackage{ amstext,amssymb,amsfonts}
\usepackage[english]{babel}
\usepackage{graphicx}
\usepackage{delarray}
 \usepackage{mathptmx}

 \DeclareFontFamily{U}{mathx}{\hyphenchar\font45}
\DeclareFontShape{U}{mathx}{m}{n}{
      <5> <6> <7> <8> <9> <10>
      <10.95> <12> <14.4> <17.28> <20.74> <24.88>
      mathx10
      }{}
\DeclareSymbolFont{mathx}{U}{mathx}{m}{n}
\DeclareFontSubstitution{U}{mathx}{m}{n}
\DeclareMathAccent{\widecheck}{0}{mathx}{"71}
\DeclareMathAccent{\wideparen}{0}{mathx}{"75}

\newtheorem{theorem}{THEOREM}
\newtheorem{proposition}[theorem]{PROPOSITION}
\newtheorem{corollary}[theorem]{COROLLARY}

   \newcommand{\R}{\mathbb R}    \newcommand{\T}{\mathbb T} \newcommand{\Co}{\mathbb C}   \newcommand{\tn}{T_N} \newcommand{\Z}{\mathbb Z}
   \newcommand{\tnp}{T_N^p}

\title{On large sieve inequalities involving $p$th powers of trigonometric polynomials}
\author{{\large S. NORVIDAS} }
\date{{\footnotesize Vilnius University Institute of Mathematics and Informatics,
 Akademijos 4,  Vilnius  08663, Lithuania\\
 {\rm{e-mail: norvidas{@}gmail.com}}}}
\begin{document}
\maketitle
\ \ \ \ \ {{ {\bf Abstract}}}.\quad In this paper, we extend the large sieve type estimates to sums involving $p$th powers of trigonometric polynomials. An approach to such estimates that does not rely  on the  usual $L^2$-technique is given.  Our method is based on  comparing the norm and the spectral radius of  convolution operators on a normed space of trigonometric polynomials.

\large
{\begin{center}
\section{ Introduction }
\end{center}}
We denote by $\tn$   the set of trigonometric polynomials of degree at  most $N$
\[
s(x)=\sum_{k=-N}^N c_ke^{ikx}
\]
with a positive integer $N$ and $c_k\in\Co$, $k=-N,\dots, N$. Suppose that   $x_1<x_2<\dots<x_r$, $r>1$,   is any sequence in $\T=(-\pi, \pi]$ such that
 \begin{equation}
\min\Bigl\{x_{j+1}-x_j, j=1,\dots, r-1; \ 2\pi-(x_r-x_1)\Bigr\}=\delta>0.
\end{equation}
The  usual large sieve inequality states that
  \begin{equation}
\sum_{j=1}^r|s(x_j)|^2\le \Bigl(\frac{N}{2\pi}+\frac1{\delta}\Bigr)\int_{\T}|s(x)|^2\,dx.
 \end{equation}
See  Selberg \cite[p. 221]{8}, but note the different notation.   Gallagher \cite{4} has given a simple derivation of the large sieve inequalities. It turned out that the method of \cite{4}  can be applied to  $L^p(\T)$-norm.   More precisely, in \cite[p. 96]{3} was proved that  if $s\in\tn$ and $2\le p<\infty$, then
\begin{equation}
\sum_{j=1}^r|s(x_j)|^p\le \Delta_pp^{1/2} \max\biggl(N, \, \frac{4\pi}{\delta}\biggr)\Bigl(\sum_{k=-N}^N|c_k|^{q}\Bigr)^{p/q},
\end{equation}
where $\Delta_p$ is an absolute constant and $1/p+1/q=1$. Next, in \cite[p. 533]{6} the inequality (3) was extended for all $0<p<\infty$ as follows. Let $\Psi$ be  a convex, nonnegative, and nondecreasing function  in $[0,\infty)$. Then for any $s\in T_N$,
  \begin{equation}
\sum_{j=1}^r\Psi(|s(x_j)|^p)\le \Bigl(\frac{N}{\pi}+\frac1{\delta}\Bigr)\int_{\T}\Psi\Bigl(|s(x)|^p(p+1)e/2\Bigr)\,dx,
 \end{equation}
whenever $-\pi<x_1< x_2<\dots <x_r\le \pi$ and $\delta$ is defined by (1). If $\Psi(t)=t$, $0\le t<\infty$, then (4) implies
\begin{equation}
\sum_{j=1}^r|s(x_j)|^p\le \Bigl(\frac{N}{\pi}+\frac1{\delta}\Bigr)\frac{(p+1)e}{2}\|s\|^p_{L^p(\T)}.
 \end{equation}
 In \cite[p. 164]{5}  the estimate (4) was extended to the case of generalized trigonometric polynomials. In particular, for usual trigonometric polynomials , the inequality (5) was improved as follows:
\begin{equation}
 \sum_{j=1}^r|s(x_j)|^p\le \Bigl(\frac{N+1}{2\pi}+\frac1{\delta}\Bigr)\frac{(p+1)e}{2}\|s\|^p_{L^p(\T)}.
 \end{equation}

Note that  inequalities (2)-(6) are also called   forward Marcinkiewicz-Zygmund inequalities (see, e.g., \cite{2}).

In this paper,   we will develop  an approach to inequalities of the type (3) and (5)-(6)  that not use the usual $L^2$-technique. Our approach is based on the  spectral theory of convolution operators on  $\tn$.

 The main result is given in the following theorem. Note that in the sequel, $[x]$ denotes the integer part of a positive number $x$. Also we use the notation $\Gamma(\cdot)$ for the standard gamma function.
\begin{theorem}.
Let $\{x_j\}_{j=1}^r$ be a sequence in $\T$  that satisfies (1). If  $s\in\tn$ and $1\le p<\infty$, then
\begin{equation}
\sum_{j=1}^r|s(x_j)|^p\le \frac{pN\sigma(\delta; N)}{2\sqrt{\pi}}\cdot \frac{ \Gamma\Bigl(p/2\Bigr)}{\Gamma\Bigl(p/2+1/2\Bigr)}\|s\|^p_{L^p(\T)}
\end{equation}
with
\begin{gather}
\sigma(\delta; N)=\begin{array}\{{rl}.
\frac{\pi}{N\delta} ,\ \ \  \ if  &  \frac{\pi}{N\delta}\in \Z, \\
1+\biggl[\frac{\pi}{N\delta}\biggr], & \mbox{otherwise.}
\end{array}
\end{gather}\
 \end{theorem}
If $p$ is a positive integer, then the quantity
\[
 \frac{ \Gamma\Bigl(p/2\Bigr)}{\Gamma\Bigl(p/2+1/2\Bigr)}
 \]
 can be calculated directly by using the relations between $ \Gamma\Bigl(p/2\Bigr)$ and $\Gamma\Bigl(p/2+1/2\Bigr)$.
 \begin{corollary}.
 Let $ l$ be a positive integer. Then under the conditions of Theorem 1 it follows that:

 (i) \ if $p=2l$, then
 \begin{equation}
\sum_{j=1}^r|s(x_j)|^p\le \frac{pN\sigma(\delta; N)}{2\sqrt{\pi}}\cdot \frac{ \Gamma\Bigl(p/2\Bigr)}{\Gamma\Bigl(p/2+1/2\Bigr)}\|s\|^p_{L^p(\T)}=\frac{pN\sigma(\delta; N)\cdot 2^{l-1}(l-1)!}{\pi\Bigl(1\cdot 3\cdot 5\cdot \cdots \cdot (2l-1)\Bigr)}\|s\|^p_{L^p(\T)};
\end{equation}
 (ii) \ \ if $p=2l+1$, then
 \begin{equation}
\sum_{j=1}^r|s(x_j)|^p\le \frac{pN\sigma(\delta; N)}{2\sqrt{\pi}}\cdot \frac{ \Gamma\Bigl(p/2\Bigr)}{\Gamma\Bigl(p/2+1/2\Bigr)}\|s\|^p_{L^p(\T)}=\frac{pN\sigma(\delta; N)\cdot \Bigl(1\cdot 3\cdot 5\cdot \cdots \cdot (2l-1)\Bigr)}{2^{l+1}\cdot l!}\|s\|^p_{L^p(\T)}.
\end{equation}
 \end{corollary}
\begin{corollary}.
Assume that $s\in\tn$,  $\|s\|_{L^p(\T)}=1$, $p\ge 1$,  and  $\{x_j\}_{j=1}^r$  satisfies (1).

(i) \  If $\frac{\pi}{N\delta}\in \Z$, then
\begin{equation}
\sum_{j=1}^r|s(x_j)|^p<\frac{ p+1}{\delta}.
\end{equation}
(ii) If $\frac{\pi}{N\delta}\not\in \Z$, then
\begin{equation}
\sum_{j=1}^r|s(x_j)|^p< (p+1)\biggl(\frac{N}{\pi}+\frac1{\delta}\biggr).
\end{equation}
\end{corollary}

\section{  Proofs}
Let $M(\T)$ be the Banach algebra of finite complex-valued regular Borel measures on $\T$. The norm in $M(\T)$ is given by the total variation $\|\mu\|$ of $\mu\in M(\T)$. Therefore, the usual Banach space $L^1(\T)$ can be identified with the closed ideal in $M(\T)$ of all measures which are absolutely continuous with respect to the Lebesgue measure $dt$ on $\T$.

Given $\mu\in M(\T)$ and $f\in L^1(\T)$, we define the Fourier transform of $\mu$ and $f$  by
\[
\widehat{\mu}(x)=\int_{\T}e^{-ixt}\,d\mu(t) \quad{\text{ and}}\quad \widehat{f}(x)=\int_{\T}e^{-ixt}f(t)\,dt,
\]
respectively.
For each $u\in L^r(\T)$, $1\le r\le\infty$, on $T_N$ is well defined  the convolution operator
\begin{gather}
A_u(s)(x)=s\ast u(x)=\int_{\T}s(x-t) u(t)\,dt= \int_{\T}\Bigl(\sum_{k=-N}^N c_ke^{ik(x-t)}\Bigr)u(t)\,dt\nonumber\\
=\sum_{k=-N}^N c_k\widehat{u}(k)e^{ikx},
\end{gather}
for each $x\in \T$, where $s\in\tn$. Note that in this definition and also below we assume  that $s$ is a periodic function on the real line with the period equal to  $2\pi$. Below,  the notation $T^p_N$, $1\le p\le\infty$, means that $T_N$ is equipped with the usual $L^p(\T)$ norm.

Suppose that $S_1$  and $S_2$ are two  measurable subsets of $\T$, $\mu_1$ and $\mu_2$ are two non-negative finite measures on $\T$, and $ F:\T^2\to\R$ is a measurable function.  For $1\le p<\infty$, Minkowski's integral inequality \cite[p. 37]{9} states that
\begin{gather}
\biggl[\int_{S_2}\biggl|\int_{S_1}F(x,y)\,d\mu_1(x)\biggr|^p\,d\mu_2(y)\biggr]^{1/p}\nonumber\\
\le\int_{S_1}\biggl(\int_{S_2}|F(x,y)|^p\,d\mu_2(y)\biggr)^{1/p}\,d\mu_1(x).
\end{gather}

\begin{proposition}.
Let $u\in L^q(\T)$, $1\le q\le \infty$, $\|u\|_{L^q(\T)}\neq 0$. Assume that $u$ is continuous,  non-negative and even on $\T$. If
\begin{equation}
{\text{supp}} \ u\subset\Bigl[-\frac{\pi}{2N}\,, \ \frac{\pi}{2N}\Bigr],
\end{equation}
then there exists a trigonometric sum
\begin{equation}
p_u(x)=\sum_{m=-N+1}^{N}\tau_m e^{-i\pi m x/N}
\end{equation}
such that
\begin{equation}
(-1)^m\tau_m> 0,
\end{equation}
$m=-N+1,\dots ,N$, and
\begin{equation}
p_u(n)=\frac1{\widehat{u}(n)}
\end{equation}
for all $n=-N,\dots,N$.
\end{proposition}
PROOF.\quad We start by examining in more details the Fourier transform of $u$. Under the assumptions on $u$, we see that
\begin{equation}
\widehat{u}(x)=\int_{\T}u(t)e^{-ixt}\,dt=2\int_{0}^{\pi/2N}u(t)\cos xt\,dt>0
\end{equation}
for all $x\in [-N, N]$. From this it follows that
\begin{equation}
(\widehat{u})'(x)=-2\int_{0}^{\pi/2N}t u(t)\sin xt\,dt<0
\end{equation}
and
\begin{equation}
(\widehat{u})''(x)=-2\int_{0}^{\pi/2N}t^2u(t)\cos xt\,dt<0
\end{equation}
for all $x\in [0, N]$.

Let $v=1/\widehat{u}$. By (19),  the function $v$  is  well defined and positive on $[-N, N]$. Moreover, we conclude from (19)-(21) that $v$ is an even function and
\begin{equation}
v'(x)>0 \quad {\text{and}}\quad v''(x)>0
\end{equation}
for all $x\in [0, N]$. Therefore,   $v$ is increasing and   convex on $[0,N]$, in particular $v$ is of bounded variation  on $[0, N]$. In particular, this means that $v$  is a function of bounded variation on $\T$. The following is well known: If $f$ is an  $2\pi$-periodic
continuously differentiable even function on $\T$ such that $f$  is of bounded variation on $\T$, then the Fourier series of $f$  converges absolutely (see, eg. \cite[p. 241]{10}).  Thus,
\begin{equation}
v(x)=\sum_{k\in\Z}a_ke^{ikx\pi/N}=\sum_{k\in\Z}a_k\cos\Bigl(\frac{\pi}{N}kx\Bigr)
\end{equation}
with
\begin{equation}
\sum_{k\in\Z}|a_k|<\infty,
\end{equation}
where
\begin{equation}
a_k=\int_{-N}^{N}v(t) e^{-ikt \pi/N}\,dt=2\int_{0}^{N}v(t)\cos\Bigl(\frac{\pi}{N}kt\Bigr)\,dt.
\end{equation}
We claim that
\begin{equation}
(-1)^ka_k>0
\end{equation}
for all $k\in\Z$. Combining (19) with (25), we see that  $a_0>0$. Let $k\ge 1$. Then using
integration by parts, we conclude from (15) that
\begin{equation}
a_k=-\frac{2N}{\pi k}\int_{0}^{N}v'(t)\sin\Bigl(\frac{\pi}{N}kt\Bigr)\,dt=-\frac{2N}{\pi k}\ \sum_{j=0}^{k-1}I_j,
\end{equation}
where
\begin{equation}
I_j=\int_{E_j}v'(t)\sin\Bigl(\frac{\pi}{N}kt\Bigr)\,dt
\end{equation}
and $E_j=[Nj/k, N(j+1)/k]$. Note that the length of $E_j$  is   exactly  half length of period for $\sin\Bigl(\pi kx/N\Bigr)$. Combining this with (22), we see that
\begin{equation}
(-1)^j I_j>0 \quad {\text{and}}\quad |I_j|<|I_{j+1}|
\end{equation}
for all $j=0,\dots,k-1$.  From this,  it is easily seen  that
\[
(-1)^{k-1}\sum_{j=0}^{k-1} I_j>0.
\]
In light of (27) this  proving the claim (26).

For $m=-N+1,\dots, N$, let $\tau_m$ be defined by
\begin{equation}
\tau_m=\sum_{j\in\Z}a_{m+2jN}.
\end{equation}
From (24), we see that the series in (30) converges absolutely. Next, from (26) it follows that,   for each $m$, all  terms of the sequence $\{a_{m+2jN}\}_{j\in\Z}$ have the same sign. In particular, (26) shows that $(-1)^ma_{m+2jN}>0$. Therefore, the trigonometric sum  (16) is well defined and satisfies (17).

Finally, for any $n\in\{-N,\dots, N\}$, combining (23) with (30), we get
\begin{gather}
\frac1{\widehat{u}(n)}=v(n)=\sum_{k\in\Z}a_ke^{ik\pi n/N}=\sum_{m=-N+1}^N\biggl(\sum_{j\in\Z}a_{m+2jN}e^{i(m+2jN)\pi n/N}\Biggr)
 \nonumber\\
=\sum_{m=-N+1}^N\biggl(\sum_{j\in\Z}a_{m+2jN}e^{im\pi n/N}\biggr)=\sum_{m=-N+1}^N\biggl(\sum_{j\in\Z}a_{m+2jN}\biggr)e^{im\pi n/N}\biggr)\nonumber\\=\sum_{m=-N+1}^N\tau_me^{imn\pi/N}=p_u(-n).\nonumber
\end{gather}
As $u$ and $\hat{u}$ are even functions we have $p_u(-n)=p(n)$ for all $n$.
Proposition 4 is proved.

We will denote by $\delta_a$ the usual Dirac measure supported on $a\in\T$.

\begin{proposition}.
Let  positive numbers $p$ and $q$ satisfy $1/p+1/q=1$. Under the conditions of Proposition 4 on  $u\in L^q(\T)$ it follows that the operator (13) possesses on $\tnp$ a bounded inverse of the type
\begin{equation}
A^{-1}_u(s)(x)=\int_{\T}s(x-t)\,d\mu(t),
\end{equation}
where
\begin{equation}
\mu=\sum_{m=-N+1}^{N}\tau_m \delta_{\pi m/N}
\end{equation}
and  $\{\tau_m\}_{-N+1}^{N}$ are defined by (30), using (23) for $v=1/\hat{u}$.  Furthermore,
\begin{equation}
\|A^{-1}_u\|_{\tnp}=\|\mu\|_{M(\T)}=\sum_{m=-N+1}^N|\tau_m|=\frac1{\widehat{u}(N)}.
\end{equation}
\end{proposition}
PROOF.\quad Let $s(x)=\sum_{k=-N}^Nc_ke^{ikx}\in \tn$. From (16) it follows that
\begin{equation}
s\ast\mu(x)=\sum _{k=-N}^Nc_k\biggl(\int_{\T}e^{ik(x-t)}\,d\mu(t)\biggr)=\sum _{k=-N}^Nc_k\widehat{\mu}(k)e^{ikx},
\end{equation}
where
\[
\widehat{\mu}(k)=\int_{\T}e^{-ikt}\,d\mu(t)=\sum_{m=-N+1}^N\tau_me^{-i\pi km/N}=p_u(k),
\]
for all $k=-N,\dots, N$. Now, taking into account (18), we conclude from from (13) that (31) defines  the inverse of the operator $A_u$.

Now,  according to  (34), we see that  the set  $\{\widehat{\mu}(k)=p_u(k): \ k=-N,\dots,N\}$ coincides with the spectrum of $A^{-1}_u$. Therefore, if $| A^{-1}_u|_{\tnp}$ denotes the spectral radius of $A^{-1}_u$, then
\[
| A^{-1}_u|_{\tnp}= \max\{|\widehat{\mu}(k)|=|p_u(k)|: \ k=-N,\dots,N\}.
\]
Let us recall that  $|H|_X\le \|H\|_X$, i.e. the spectral radius is not greater than the operator norm,  for any bounded linear operator $H$ on a normed space  $X$.
Combining this with (18), we get
\begin{gather}
\|A^{-1}_u\|_{\tnp}\ge \max_{-N\le k\le N}|\widehat{\mu}(k)|= \max_{-N\le k\le N}|p_u(k)|=\max_{-N\le k\le N}\biggl\{\frac1{|\widehat{u}(k)|}\biggr\}\nonumber\\
=\frac1{\min_{-N\le k\le N}|\widehat{u}(k)|}
\ge\frac1{\widehat{u}(N)}=p_u(N)=\sum_{M=-N+1}^N \tau_me^{i\pi mN/N}=\nonumber\\
\sum_{M=-N+1}^N (-1)^m\tau_m=\sum_{M=-N+1}^N |\tau_m|=\|\mu\|.
\end{gather}
Other hands, from Minkowski's inequality (14) it is easily to see that
\begin{gather}
\|A^{-1}_u(s)\|_{\tnp}=\biggl(\int_{\T}\biggl|\int_{\T}s(x-t)\,d\mu(t)\biggr|^p\,dx\biggr)^{1/p}\nonumber\\
\le \biggl(\int_{\T}\biggl|\int_{\T}|s(x-t)|\,d|\mu|(t)\biggr|^p\,dx\biggr)^{1/p}
\nonumber\\
\le \int_{\T}\biggl(\int_{\T}|s(x-t)|^p\,dx\biggr)^{1/p}\,d|\mu|(t)= \|s\|_{\tnp}\|\mu\|,
\end{gather}
where $|\mu|$ denote the variation of $\mu$. Thus, (35) with (36) show   (33), and the Proposition 5  is proved.

PROOF OF THEOREM 1.\quad
Assume that $u\in L^q(\T)$, $\|u\|_{L^q(\T)}=1$ and $u$ satisfies the conditions of Proposition 4.  If $s\in\tn$, then by H\"{o}lder's inequality, we get
\begin{gather}
|s(x)|^p= \biggl|\int_{\T}A^{-1}_u(s)(x-t)u(t)\,dt\biggr|^p= \biggl|\int_{-\pi/2N}^{\pi/2N}A_u^{-1}(s)(x-t)u(t)\,dt\biggr|^p\nonumber\\
\le \biggl(\int_{-\pi/2N}^{\pi/2N}\Bigr|A_u^{-1}(s)(x-t)\biggl|^p\,dt\biggr) \biggl(\int_{-\pi/2N}^{\pi/2N}|u(t)|^q\,dt\biggr)^{p/q}\nonumber\\
= \int_{-\pi/2N}^{\pi/2N}\Bigr|A^{-1}_u(s)(x-t)\biggl|^p\,dt\nonumber
\end{gather}
for each $x\in\T$. Therefore,
\begin{gather}
\sum_{j=1}^r|s(x_j)|^p\le \sum_{j=1}^r\int_{-\pi/2N}^{\pi/2N}\biggl|A_u^{-1}(s)(x_j-t)\biggr|^p\,dt\nonumber\\ =\sum_{j=1}^r\int_{x_j-\pi/2N}^{x_j+\pi/2N}\biggl|A_u^{-1}(s)(y)\biggr|^p\,dy.
\end{gather}
We will denote by $E_j$ the set
\[
E_j=\Bigl(x_j-\frac{\pi}{2N} , x_j+\frac{\pi}{2N}\Bigr],
\]
$ j=1,\dots, r$. Using the fact that $A_u^{-1}(s)$ is a trigonometric polynomial, i.e., a continuous and periodical function on $\R$, we conclude that
\begin{equation}
\int_{x_j-\pi/2N}^{x_j+\pi/2N}\Bigl|A_u^{-1}(s)(y)\Bigr|^p\,dy  = \int_{E_j}\Bigl|A_u^{-1}(s)(y)\Bigr|^p\,dy.
\end{equation}
 Let $x\in\R$ and assume that $x\in E_j$, for $j=i, i+1,\dots, i+k$ with some  non-negative integer $k$.  Then we claim that:

 (i) \ if $\pi/N\delta\in\Z$, then
 \begin{equation}
k+1\le \frac{\pi}{N\delta};
\end{equation}
 (ii) \ if $\pi/N\delta\not\in\Z$, then
\begin{equation}
k\le  \Bigl[\frac{\pi}{N\delta}\Bigr],
\end{equation}
where $[\cdot]$ is the usual integer part of a  real number.

Indeed, since $x\in E_i\cap E_{i+k}$, it follows that $x_i+\pi/2N> x_{i+k}- \pi/2N$. Then
\begin{equation}
x_{i+k}-x_i< \frac{\pi}N.
\end{equation}
Other hands, by (1)  we get
\begin{equation}
x_{i+k}-x_i\ge k\delta.
\end{equation}
Combining (42) with (41), we see that
\begin{equation}
k< \frac{\pi}{N\delta}.
\end{equation}
Now note  that $k$ is an integer. Hence if $\pi/N\delta\in\Z$, then (43) implies (39). For $\pi/N\delta\not\in\Z$, we get (40), which yields our  claim.

Thus, each $x\in\T$ can  belong to at most  $\sigma_p$ intervals $E_j$, $j=1,\dots, p$,  where $\sigma_p$ was defined by (8).

Now combining (35), (36) with (37), (38), (31) and (8), we get
 \begin{gather}
\sum_{j=1}^r|s(x_j)|^p\le\sum_{j=1}^r\int_{x_j-\pi/2N}^{x_j+\pi/2N}\biggl|A_u^{-1}(s)(y)\biggr|^p\,dy\le \sigma(\delta; N)\int_{\T}\Bigl|A^{-1}_u(s)(y)\Bigr|^p\,dy=\nonumber\\
= \sigma(\delta; N) \|A^{-1}_u(s)\|^p_{L^p(\T)}\le \sigma(\delta; N)\|A_u^{-1}\|^p_{T^p_N }\|s\|^p_{L^p(\T)}=
 \frac{\sigma_p}{\bigl(\hat{u}(N)\bigr)^p}\|s\|^p_{L^p(\T)}.
\end{gather}
Next, the estimate (44) can be improved to
\begin{equation}
\sum_{j=1}^r|s(x_j)|^p\le \frac{\sigma(\delta; N)}{\sup_{u}\Bigl(\hat{u}(N)\Bigr)^p}\ \|s\|^p_{L^p(\T)},
\end{equation}
where the supremum extends over all admissible $u$ as described in the statement of Proposition 4. We claim that for such an $u$ we have
\begin{equation}
\sup_{u}(\widehat{u}(N))^p=\|\cos Nx\|^p_{L^p[-\pi/2N, \pi/2N]}=\frac{\|\cos t\|^p_{L^p[-\pi/2, \pi/2]}}{N}.
\end{equation}
Indeed, H\"{o}lder's inequality implies that
\begin{gather}
|\widehat{u}(N)|^p= \Bigl|\int_{\T}u(x)\cos Nx\  dx\Bigr|^p= \Bigl|\int_{-\pi/2N}^{\pi/2N}u(x)\cos Nx\  dx\Bigr|^p\nonumber\\
\le \int_{-\pi/2N}^{\pi/2N}|\cos Nx|^p\  dx \cdot \|u\|_{L^q[-\pi/2N, \pi/2N]}^p=\int_{-\pi/2N}^{\pi/2N}|\cos Nx|^p\,dx\nonumber\\
=\frac1{N}\int_{-\pi/2}^{\pi/2}|\cos t|^p\, dt=\frac{\|\cos t\|^p_{L^p[-\pi/2, \pi/2]}}{N}.
\end{gather}
 Moreover, since we used H\"{o}lder's inequality, it follows that the estimate (47) is exact and the equality is attained if
\[
u(x)= \theta\cos ^{p-1} Nx \cdot\chi_{[-\pi/2N, \pi/2N]}(x),
\]
 where $\chi_{[-\pi/2N, \pi/2N]}$ is the indicator function of the interval $[-\pi/2N, \pi/2N]$  and $\theta \in\Co$ is such that   $\|u\|_{L^q[-\pi/2N, \pi/2N]}=1$. Therefore, our claim (46) is proved. 
 
 Combining (45) with (46), we get  
 \begin{equation}
 \sum_{j=1}^r|s(x_j)|^p\le\frac{\sigma(\delta; N)\cdot N}{\|\cos t\|^p_{L^p[-\pi/2, \pi/2]}}\|s\|^p_{L^p(\T)}.
 \end{equation}
 Next,
  \begin{equation}
\|\cos t\|^p_{L^p[-\pi/2, \pi/2]}=\int_{-\pi/2}^{\pi/2} \cos^p t\, dt=B\Bigl(\frac12; \frac{p+1}{2}\Bigr),
 \end{equation}
 (see, e.g., \cite[p. 142]{7}), where $B$ is Euler's beta function defined by 
 \[
 B(a;b)=2\int_0^{\pi/2}\sin^{2a-1}\theta\cos^{2b-1}\theta\, d\theta
 \]
 for $\Re a$, $\Re b>0$.  Applying the following connection between the beta and the usual gamma function  (\cite[p. 142]{7}
 \[
 B(a,b)=\frac{\Gamma(a)\Gamma(b)}{\Gamma(a+b)}, 
 \]
 we conclude from (49) that 
  \begin{equation}
\|\cos t\|^p_{L^p[-\pi/2, \pi/2]}=\frac{\Gamma\Bigl(\frac12\Bigr)\Gamma\Bigl(\frac{p}2+\frac12\Bigr)}{\Gamma\Bigl(\frac{p}2 +1\Bigr)}
 \end{equation}
Since $\Gamma\Bigl(p/2+1\Bigr)=\Gamma\Bigl(p/2\Bigr)\cdot p/2$ and $\Gamma\Bigl(1/2\Bigr)=\sqrt{\pi}$ (see, e.g., \cite[p.p.  137-138]{7}), it follows from (49) and (50) that 
\[
\|\cos t\|^p_{L^p[-\pi/2, \pi/2]}=\frac{2\sqrt{\pi}\cdot \Gamma\Bigl(\frac{p}2+\frac12\Bigr)}{\Gamma\Bigl(\frac{p}2\Bigr)}.
 \] 
Substituting this into (48), we obtain (7). Theorem 1 is proved.

POOF OF COROLLARY 2. \quad It is known that for a nonnegative integer $n$,  
\[
\Gamma(n+1)=n! \qquad {\text{and}} \qquad \Gamma\Bigl(n+\frac12\Bigr)=\sqrt{\pi}\cdot\frac{1\cdot 3\cdot 5\cdot \cdots \cdot(2-1)}{2^n}
\]
 (see \cite[p.  139]{7}). Using this, the representations (9) and (10) can be verified by straightforward calculation.

\end{document}